\documentclass[10pt,reqno]{amsart}
\usepackage{amsmath,amssymb,graphicx}
\usepackage[usenames,dvipsnames]{color}

\numberwithin{equation}{section}
\numberwithin{figure}{section}
\numberwithin{table}{section}
\newtheorem{theorem}{Theorem}[section]

\newtheorem{corr}{Corollary}[section]

\newtheorem{lema}{Lemma}[section]

\newtheorem{rem}{Remark}[section]

\newcommand{\ignore}[1]{}{}

\usepackage{amssymb,color}
\definecolor{c20}{rgb}{0.,0.7,0.}
\definecolor{c30}{rgb}{0.,0.,1.}
\definecolor{c40}{rgb}{1,0.1,0.7}
\definecolor{c50}{rgb}{1,0,0}

\def\aa#1{\textcolor{black}{#1}} 
\def\aa#1{#1}

\def\pE#1{\textcolor{c20}{#1}}
\def\pE#1{#1}

\def\cE#1{#1}

\def\tE#1{\textcolor{c20}{#1}}
\def\tE#1{#1}

\newcommand{\BQN}{\begin{eqnarray}}
\newcommand{\EQN}{\end{eqnarray}}
\newcommand{\BQNY}{\begin{eqnarray*}}
\newcommand{\EQNY}{\end{eqnarray*}}

\newcommand{\pk}[1]{\mathbb{P}\left\{#1\right\} }

\newcommand{\R}{\mathbb{R}}

\newcommand{\inr}{\in \R}

\newcommand{\prooftheo}[1]{ \textsc{Proof of Theorem} \ref{#1} }
\newcommand{\proofprop}[1]{\textsc{Proof of Proposition} \ref{#1}}
\newcommand{\prooflem}[1]{\textsc{Proof of Lemma} \ref{#1}}
\newcommand{\proofkorr}[1]{\textsc{Proof of Corollary} \ref{#1}}
\newcommand{\QED}{\hfill $\Box$}
\newcommand{\IF}{\infty}
\newcommand{\BL}{\begin{lem}}
\newcommand{\EL}{\end{lem}}
\newcommand{\nelem}[1]{{Lemma \ref{#1}}}

\newcommand{\netheo}[1]{{Theorem \ref{#1}}}

\newcommand{\ldot}{,\ldots,}

\def\limit{\lim_{u\to \IF}}

\newcommand{\abs}[1]{\lvert #1 \rvert}

\def\ve{\varepsilon}

\newcommand{\ABs}[1]{ \biggl \lvert #1 \biggr \rvert}

\newcommand{\BT}{\begin{theorem}}
\newcommand{\ET}{\end{theorem}}
\newcommand{\BEL}{\begin{lema}}
\newcommand{\EEL}{\end{lema}}

\newcommand{\equaldis}{\stackrel{d}{=}}

\usepackage{color}
\newcommand{\tb}[1]{{\textcolor{blue}{#1}}}
\newcommand{\tc}[1]{{\textcolor{blue}{#1}}}
\def\tb#1{#1}
\def\tE#1{#1}
\def\tc#1{#1}

\begin{document}

\vspace*{1cm}

\begin{center}

\MakeUppercase{\bf Random Scaling of Gumbel Risks}

\large \normalsize

\bigskip

\centerline{\large
Krzysztof D\c{e}bicki\footnote{Mathematical Institute, University of Wroc\l aw, pl. Grunwaldzki 2/4, 50-384 Wroc\l aw, Poland},
Julia Farkas$^{2}$ and Enkelejd  Hashorva\footnote{University of Lausanne, Faculty of Business and Economics (HEC Lausanne),
Lausanne 1015, Switzerland}}

\bigskip

\today{}

\bigskip

\end{center}

\begin{quote}
\noindent
\textbf{Abstract.} In  this paper we consider the product of two positive independent risks  $Y_1$ and $Y_2$. If $Y_1$ is bounded and $Y_2$ has distribution in the Gumbel max-domain of attraction with some auxiliary function which is regularly varying at infinity,
then we show that $Y_1Y_2$ has also distribution in the Gumbel max-domain of attraction. Additionally, if both $Y_1,Y_2$ have  log-Weibullian or Weibullian tail behavior, we show that $Y_1Y_2$ has log-Weibullian or Weibullian \tE{asymptotic} tail behavior, respectively.

\bigskip

\noindent
\textit{Keywords and phrases}: Gumbel max-domain of attraction; product of random variables; log-Weibullian tail behavior;
Weibullian tail behavior; supremum of Gaussian process.
\end{quote}

\section{Introduction}

Consider $Y_1$ and $Y_2$, two positive independent random
variables  \pE{(rvs)}. If $Y_2$ is bounded, say $Y_2\le 1$ almost
surely, then $X=Y_1Y_2$ is referred to as a random contraction,
see e.g., Pakes and Navarro (2007). In such a contraction model we
expect that
 the asymptotic tail behavior of $X$ is essentially determined by that of $Y_1$.
 This intuition is confirmed in \netheo{thm1} below for the case $Y_1$ has a
 distribution with unbounded support, being further in the Gumbel max-domain of attraction, i.e.,
\BQN\label{conG}
\limit \frac{ \pk{ Y_1> u+  a(u)t}}{\pk{Y_1> u}}= \exp(- t), \quad \forall t\ge 0
\EQN
\tb{for some positive scaling function $a(\cdot)$, which is regularly varying at infinity
with index $-\tau$ for $\tau\pE{\ge -1}$, i.e., $\limit a(u x)/a(u)= x^{-\tau},x>0$.
We abbreviate \aa{(1.1)} as $Y_1 \in GMDA(a)$ and refer to, e.g., Resnick (1987),
for details on the Gumbel max-domain of attraction and regular variation.}

\BT \label{thm1}
If condition \eqref{conG} holds with $a(\cdot)$ being
regularly varying at infinity with index $-\tau$ for $\tau\pE{\ge -1}$ and \tc{$Y_2$} \pE{has distribution  with right endpoint equal to 1}, then $X=Y_1Y_2 \in GMDA(a)$.
\ET
In view of Lemma 3.2 in Arendarczyk and D\c{e}bicki (2011) a direct implication of \netheo{thm1}  is
 $$X=\sup_{t\in [0,S]} B_H(t) \in GMDA(a), \quad \text{ where } a(x)=1/x,$$
provided that  $S>0$ is a bounded risk being independent of
\tb{
a standard fractional Brownian motion $\{B_H(t), t\inr\}$}
(\tE{with} mean zero, variance function $t^{2 H}$) with Hurst \pE{index} $H\in (0,1)$.\\
A canonical example for  $Y_1 \in GMDA(a)$ is when  
\begin{eqnarray}\label{poly}
\pk{Y_1> u} \sim C_1 u^{\alpha_1} \exp(- L_1 u^{p_1}),  \quad u\to \IF,
\end{eqnarray}
where  $C_1,L_1,p_1$ are positive constants and $\alpha_1 \inr$;
\tE{note} that $f_1(u) \sim f_2(u)$ means  $\limit f_1(u)/f_2(u)=1$.\\
Clearly, if \eqref{poly} holds, then $Y_1 \in GMDA(a)$, where  $a(u)= u^{1- p_1}/L_1$. Consequently,  the assumption of \netheo{thm1} on $a(\cdot)$ holds with $\pE{\tau=p_1-1}$.

If $Y_1$ and $Y_2$ can simultaneously take large values with non-zero probability, then the asymptotic tail behavior of $X$ is known in few cases. In particular, if \tE{also $Y_2$ satisfies \eqref{poly} with some constants $\alpha_2\inr ,C_2>0, L_2>0, p_2>0$},
then  in light of Arendarczyk and D\c{e}bicki (2011)
\begin{eqnarray}\label{DE}
\pk{X> u}\sim\Bigl(\frac{2\pi p_{2}L_{2}}{p_{1}+p_{2}}\Bigr)%
^{1/2}C_{1}C_{2}A^{p_{2}/2+\alpha _{2}-\alpha _{1}}u^{\frac{2p_{2}\alpha
_{1}+2p_{1}\alpha _{2}+p_{1}p_{2}}{2(p_{1}+p_{2})}}
 \exp \left( -Bu^{\frac{p_{1}p_{2}}{%
p_{1}+p_{2}}}\right)
\end{eqnarray}
holds as $u \to \IF$, where
 \BQN\label{CA}
 A =[(p_1L_1)/(p_2 L_2)]^{1/(p_1+p_2)} \mbox{   and   } B=L_{1}A^{-p_{1}}+L_{2}A^{p_{2}}.
 \EQN
Our second result shows that the asymptotic \pE{tail} behavior of  $X$ can also be derived for
\tb{a} more general case
when the power term in the tail expansion of $Y_1$ and $Y_2$ is substituted by some regularly
varying function, see \netheo{th2} in Section 2.
\tb{We refer to, e.g., Berman (1992), Cline and Samorodnitsky (1994),
Maulik and Resnick (2004), Nadarajah (2005), 
D\c{e}bicki  and van Uitert (2006), Jessen and Mikosch (2006),  D'Auria  and Resnick (2006, 2008), Resnick (2007), Charpentier and Segers (2007, 2009), Schlueter and Fischer (2012), D\c{e}bicki et al.\ (2013),
\pE{Farkas and Hashorva (2013)}, Hashorva and Weng (2013), Tan and Hashorva (2013),  Yang and Wang (2013) for related
results and numerous motivations of investigation of tail behavior of the distribution of products of rvs.
}

\tb{
As an illustration of Theorems \ref{thm1} and \ref{th2}, we analyze:\\
$\diamond$ limiting behavior of the maximum of randomly scaled Gaussian processes,\\
$\diamond$ exact asymptotic tail behavior of the supremum of Gaussian processes with stationary increments
over a random interval with length which has
Weibullian tail behavior.}

We organize this paper as follows: Section 2 derives the
 tail \aa{asymptotics} of the product of two independent
Weibullian-type rvs. Our applications are presented
in Section 3. \aa{Proofs} of all results are relegated to
Section 4, which concludes this article.

\section{Log-Weibullian and Weibullian Risks}
We say that $Y_i,i=1,2$ has a {\it{log-Weibullian}} tail behavior (or alternatively $Y_i$ is a log-Weibullian rv), if
\BQN\label{logW}
\lim_{u\to\infty} \frac{\log(\pk{Y_i>u})}{u^{p_i}}=-L_i
\EQN
for some positive constants $p_i,L_i$. The main result in this section is \netheo{th2};
statement (a) therein shows that if \eqref{logW} holds, then $X=Y_1Y_2$ has also a log-Weibullian tail behavior.

The definition of Weibullian tail behavior is  formulated (motivated by \eqref{poly}) in terms of the following condition:
\BQN\label{giv}
  \pk{Y_i> u} \sim  g_i(u) \exp(- L_i u^{p_i}),  \quad u\to \IF
\EQN
holds for $i=1,2$, where  $g_1,g_2$ are two given
\tc{regularly varying at infinity} functions and $L_i,p_i,i=1,2$ are positive constants. We say alternatively that $Y_1$ and $Y_2$ are Weibullian-type rvs. \tc{We note that if a rv is of Weibullian-type, then it is log-Weibullian}.

 For  $g_1,g_2$ being regularly varying and ultimately monotone Hashorva and Weng (2013) shows that a similar result to
\eqref{DE} holds. In statement (b) of \netheo{th2}
we remove the assumption that $g_1$ and $g_2$
are ultimately monotone.

\newcommand{\COM}[1]{}

\BT\label{th2}
Let $Y_1,Y_2$ be two independent positive rvs, and let $L_i,p_i,i=1,2$ be positive constants.\\
(a) If $Y_i,i=1,2$ satisfy \eqref{logW} with  $p_i,L_i,i=1,2$, then
with  $B$ given in \eqref{CA}
\BQN
\lim_{u\to\infty}
\frac{\log(\pk{Y_1Y_2>u})}{u^{p_1p_2/(p_1+p_2)}}=-B.
\EQN
(b) Assume that for $i=1,2$ \tE{condition \eqref{giv} holds with $L_i,p_i,i=1,2$ and}
 $g_1,g_2$ two regularly varying functions at infinity. If $A,B$ are two constants as
in \eqref{CA}, then we have
\BQN
\pk{Y_1Y_2>u}&\sim& D
u^{\frac{p_{1}p_{2}}{2(p_{1}+p_{2})}}g_1(u/c_u)g_2(c_u)
\exp\left(-Bu^{\frac{p_1p_2}{p_1+p_2}}\right)\notag \\
& \sim & D
u^{\frac{p_{1}p_{2}}{2(p_{1}+p_{2})}} \pk{Y_1 > u/c_u} \pk{Y_2 > \tE{c_u}}
\EQN
as $u\to \IF$,  where
$$c_u=Au^{p_1/(p_1+p_2)}, \quad \text{ and } D=\Biggl(\frac{2\pi  (p_1L_1)^{\frac{p_{2}}{p_1+p_2}}  (p_2L_2)^{\frac{p_{1}}{p_1+p_2}} }{p_{1}+p_{2}}   \Biggr)^{1/2}.$$
\ET

\bigskip

\begin{rem}
Theorem \ref{th2} straightforwardly extends to the case of the product of $n$ rvs.
Namely, if $Y_i,i\le n$ are positive independent rvs with tail asymptotics given by \eqref{giv}, then
$X=\prod_{i=1}^n Y_i$  also satisfies the condition \eqref{giv} with some $g^*, L^*$ and $p^*= (\sum_{i=1}^n 1/p_i)^{-1}$.
\end{rem}

\tc{Hereafter by $h^{\leftarrow}(u):=\inf\{x:h(x)\ge u\}$ we denote the generalized inverse
of the function $h$.}

\begin{rem}

\tc{Let $(Y_{n1},Y_{n2}),n\ge 1$ be independent copies of $(Y_1,Y_2)$
and let
$F^{\leftarrow}$  and $H^{\leftarrow}$
be the generalized inverse of the distributions of $Y_1$ and $Y_1Y_2$, respectively.
}
Define next 
\BQNY
\tc{b_n:=F^{\leftarrow}(1-1/n), \quad \widetilde b_n:=H^{\leftarrow}(1-1/n), \quad n>1.}
\EQNY
\tc{Under the assumptions of \netheo{thm1} for $Y_1$ and $Y_2$, we
have that $a(b_n) \sim a(\widetilde b_n)$ and $\widetilde b_n\sim b_n$ and}
further 
\BQNY
\lim_{n \to \IF }\sup_{x\inr} \ABs{\pk{\max_{1 \le k\le n}  Y_{k2} \le a(b_n) x+ b_n} - \exp(-\exp(-x))}=0.
\EQNY
In view of \netheo{thm1} it follows that $Y_1Y_2$ is in Gumbel \aa{max-domain of attraction}. Furthermore,
since
\tc{$H^{\leftarrow}(1-1/x)$}
is a slowly varying function at infinity, we obtain
\BQNY
\lim_{n \to \IF }\sup_{x\inr} \ABs{\pk{\max_{1 \le k\le n}  Y_{k1}Y_{k2} \le a(b_n) x+ \widetilde b_n} - \exp(-\exp(-x))}=0.
\EQNY
\COM{ii) Under additional conditions on the existence of the probability density functions for $Y_1$ and $Y_2$, it is possible
to derive the asymptotic behavior of the probability density function of the product $Y_1Y_2$.\\
}
\COM{ii) Our results can be shown to hold also for the case that $Y_1,Y_2$ are both non-negative, or $Y_1$ is non-negative and $Y_2$ is
a real-valued rv. Furthermore, the assumption of independence of $Y_1$ and $Y_2$ can
be removed considering some tractable classes such as FGM-type dependence, see Hashorva and Weng (2013).
}
\label{RemA}
\end{rem}

\section{Applications}
In this section we \pE{present} two applications of \netheo{thm1} and \netheo{th2}.
The first
one focuses on
the maximum of randomly scaled Gaussian processes.
\COM{extremes of establishes the Gumbel MDA of infinite
sum of inflated products of independent Gaussian rvs.
Utilizing the results of both \netheo{thm1} and \netheo{th2} in
our second application we are able to deal with the tail
asymptotic behavior of elliptical chaos. Continuing, we present
three other applications of both aforementioned theorems.
The most significant application
presented in this paper  is the fifth (last) one,
}
The second one,
\tb{which combines}
\netheo{th2}  with \pE{an interesting} finding of Arendarczyk and D\c{e}bicki (2011), \aa{derives} the asymptotic
behavior of the tail distribution of supremum
of Gaussian processes with stationary increments over Weibullian and log-Weibullian random intervals.
\COM{
\bigskip
\subsection{Infinite sum of deflated Gaussian products}
Consider $R_i,X_i,Y_i,i\ge 1$ mutually independent rvs  such that
$X_i,Y_i$ have $N(0,1)$ distribution for any index $i$ and $R_i\in (0,1]$ almost surely. Hashorva et al.\ (2012)
derived the tail asymptotics of
$$ S_\IF:= \sum_{i=1}^\IF \lambda_i R_i X_iY_i$$
assuming that $\lambda_i$'s are non-negative weights with $\sum_{i=1}^\IF \lambda_i^2 \in (0,\IF)$ and for $R_i$'s the regular variation at the right endpoint 1 was further assumed. We relax the latter assumption on \tE{$R_i$'s} considering more special case of $\lambda$'s:
\BQN\label{Li}
1=\lambda_1 = \cdots =\lambda_m, \text{ with } m\ge 1, \text{ and  }\lambda_{i}\in [0, 1), i\ge m+1.
\EQN
\tb{
The price to pay for this extension is that instead of the asymptotic tail behavior of $S_\IF$,
we claim its belonging to the GMDA.
}

\BS \label{infsum}
\tb{
Assume that $\lambda_1,\lambda_2,...$ satisfy \eqref{Li} for some $m\ge1$.
}
Let  $X_i,Y_i,i\ge 1$ be independent $N(0,1)$ rvs  being further independent of $R_i,i\ge 1$. Suppose that
$\sum_{i=m+1}^\IF \lambda_i^2 \in (0,\IF)$ and $R_i\in (0,1],i\ge 1$ almost surely. If further $R_1=R_2 =\cdots =R_m$ \tE{and $R_i,i\le m$ are all independent}, then
$ S_\IF \in GMDA(a)$  \tb{with $a(x)=1$ for $x>0$}.
\ES

\bigskip
\subsection{Gumbel MDA for elliptical chaos}  For notation simplicity we consider only two dimensional setup.
\tb{Let $(X_1,X_2)$} be a bivariate random vector with stochastic representation
\BQN \label{stocR} (X_1,X_2) \equaldis \Bigl(R O_1, R (\rho O_1+ \sqrt{1- \rho^2} O_2)\Bigr),
\EQN
where $\rho \in (-1,1)$ and $\equaldis$ means equality of distributions.  For the random radius $R>0$ we shall assume
that  it has some distribution  $F$ with right \tb{endpoint \pE{equal} infinity}, such that $F \in GMDA(a)$.
Further we suppose that $(O_1,O_2)$ is independent of $R$ and $O_1^2+O_2^2=1$ almost surely.
In the particular case that $(O_1,O_2)$ is uniformly distributed on the unit circle $\{(x,y)\inr^2: x^2+y^2=1\}$
\tb{(and thus $O_1$ is symmetric about 0 with $O_1^2$ having  beta distribution with parameters (1/2,1/2)),
random vector
$(X_1,X_2)$ is elliptically symmetric; for more details see Cambanis et al.\ (1981).}
Note that in the elliptical setup, if $R^2$ is chi-square distributed, then $(X_1,X_2)$ is a bivariate Gaussian random vector with $N(0,1)$ components and correlation $\rho$.

Next, for some positive homogeneous function $h(x_1,x_2)$ of order $\lambda >0$ i.e.,
$$h(cx_1,c x_2)= c^\lambda h(x_1,x_2), \quad  \text{ for any  } c>0,x_1,x_2\inr$$
 we define the rv
$$ Z= h(X_1,X_2),$$
which is referred to as elliptical  chaos.
\tb{ Some particular examples of $h$ satisfying the above condition are:}  \\
i) $h_1(x_1,x_2)= (x_1 x_2)^\lambda; $\\
ii)  $h_2(x)= \abs{x_1}^\lambda + \abs{x_2}^\lambda$;\\
iii)  $h_3(x)= \min(\abs{x_1}^\lambda, \abs{x_2}^\lambda)$.\\
In the special case that $R^2$ is chi-square distributed, \tb{$Z$} is referred to as the {\it{Gaussian chaos}}; see Lata{\l}a (1999, 2006).\\
In order to determine the tail asymptotic behavior of $Z$, specific assumptions on $R$ and $(O_1,O_2)$ are needed.
Using the stochastic representation \eqref{stocR} and the homogeneity property of $h$, we may write
$$ Z \equaldis R^\lambda h(O_1, \rho O_1+ \sqrt{1-\rho^2} O_2)=: R^\lambda U.$$
Consequently, if the scaling function $a(\cdot)$ is regularly varying at infinity, then
we have $R^\lambda \in GMDA(a^*)$ with
$$a^*(x)=  \frac{1}{\lambda} x^{1/\lambda -1} a(x^{1/\lambda}), \quad x>0.$$
Since also $a^*$ is a regularly varying function at infinity, if further $h$ is such that $U$ is a positive bounded
rv, then \netheo{thm1} implies
\BQN \label{ZZ}
Z \in GMDA(a^*).
\EQN
If $(X_1,X_2)$ is a bivariate Gaussian random vector, then $R$ is in Gumbel MDA with scaling function $a(x)=1/x$, and consequently  \eqref{ZZ} holds for this case.\\
Clearly, if $R$ and $U$ are two Weibullian-type rvs, the exact tail asymptotics of $Z$ follows by \netheo{th2}.

}

\subsection{Limit law of the maximum of deflated Gaussian processes}
 This application is motivated by a key finding of Kabluchko (2011). Instead of Gaussian processes treated therein, we consider here
 deflated Gaussian processes. Let therefore
$\Gamma(\cdot,\cdot)$ be a negative definite kernel in $\R^2$ and define a Brown-Resnick stochastic process with Gaussian points
as
 \BQN \label{def:zeta}
\eta_{BR}(t) =  \max_{i\ge 1} \Bigl(\Upsilon_i+ Z_i(t) -\sigma^2(t)/2  \Bigr),  \quad t \in \R  ,
  \EQN
where  $\{Z_i(t), t\inr\}, i\ge 1$
\tb{are mutually independent} centered Gaussian processes with
incremental variance function $Var(Z_i(t_1)- Z_i(t_2))=\Gamma(t_1,t_2), i\ge 1$
and variance function $\sigma^2(\cdot)>0$, being further
independent of the Poisson point process $ \Xi =\sum_{i=1}^\IF \ve_{\Upsilon_i}$  with intensity measure $\exp(-x)\, dx$, see for more details Kabluchko (2011).

In the following, for scaling the Gaussian process, we shall use a
generic positive rv $S$, which has either a
distribution with right endpoint 1, or it has a
Weibullian tail behavior satisfying \eqref{giv} with some
\tc{$p,L$ and $g$} being regularly varying at infinity. Our next
result generalizes Theorem 5.1 in Hashorva (2013).

\BT \label{th:last}
Let $\{X_{ni}(t), t\in \R\}, 1\le i\le n, n\ge 1$ be independent Gaussian processes with mean-zero,
unit variance function and
correlation function $\rho_n(s,t), s,t\inr$. Let $S_{ni},i,n\ge 1$ be
independent copies of $S$, and let
\tc{$H^{\leftarrow}$}
be the generalized inverse of the distribution
 $H$ of $S X_{11}(1)$. Assume that $S,S_{ni}, X_{ni}(t),t\in \R$ are
\tb{mutually}
independent for any $i=1 \ldot n$.
\tc{For
$  d_n= H^{\leftarrow}(1- 1/n)$ set  $c_n=1/d_n$ if $S$ is bounded, and set
\tc{$c_n= d_n (2+p)/(2p  \log  n)$}
otherwise.} If further
\BQN \label{rho:fish}
\lim_{n\to \IF} \frac{ 2 d_n}{c_n} 
\Bigl(1- \rho_n(t_1,t_2)\Bigr) &= &\Gamma(t_1,t_2)\in (0,\IF), \quad t_1\not=t_2 \inr,
\EQN
then, as  $n\to \IF$
\BQN\label{claim:M}
c_n \Bigl( \max_{1 \le i \le n} S_{ni}X_{ni}(t) -   d_n \Bigr)&\Longrightarrow &  \eta_{BR} (t), \quad t\inr,
\EQN
where $\Longrightarrow $ means the weak convergence of the finite-dimensional distributions. Furthermore, $  d_n=(1+o(1))\sqrt{2  \log  n}$ if
$S$ is bounded and
\tc{$  d_n= (1+o(1))(( \log  n)/B)^{(2+p)/(2p)}$}
otherwise with $B$ as in \eqref{CA}.
\ET

\COM{
\bigskip

\subsection{Weak tail dependence coefficient of scaled elliptically symmetric random vectors}
When a bivariate random vector $(X_1,X_2)$ with
\tb{some}
marginal distributions exhibits asymptotic independence,
a property well-established in multivariate extreme value theory (see e.g., Resnick (1997), Falk et al.\ (2010) or Grothe (2012)), then a further quantification of the strength of the tail dependence is achieved by calculating the residual
tail dependence index (see e.g., Falk et al.\ (2010)) or the weak dependence coefficient
$\bar \chi_{(X_1,X_2)}$ introduced in Coles et al. \ (1999),
defined as the  limit  (if its exits)
$$ \bar \chi_{(X_1,X_2)}= \lim_{u\to \IF}  \frac{2 \log \pk{X_1> u}}{ \log \pk{X_1> u, X_2>u}}.$$
Next, let $S$ be as in \tb{Section 3.3, i.e., either has}  distribution with right endpoint equal 1 or it has a Weibullian tail
with parameters $p_1,L_1$ and some regularly varying function $g_1$.\\
If $S$ is further independent of $(X_1,X_2)$, then  $(Y_1,Y_2)= (S X_1, S X_2)$ is a bivariate rv obtained by randomly scaling $X_1$ and $X_2$, which can be interpreted as a simple probabilistic model, consisting of a common random shock $S$ and an initial risk $(X_1,X_2)$.
It is of some interest to investigate the
impact of the common shock on the tail dependence of $(S X_1, S X_2)$.\\
Next, we consider a bivariate elliptically symmetric random vector $(X_1,X_2)$ with stochastic representation \BQN\label{Ellip}
(X_1,X_2)\equaldis R(O_1,\rho O_1+\sqrt{1-\rho^2}O_2), \qquad \qquad \rho \in(-1,1),
\EQN
where $R>0$ is independent of $(O_1,O_2)$ which is uniformly
distributed on the unit circle of $\R^2$; see Cambanis et al. (1981) for
the basic properties of elliptically symmetric random vectors.

Now, in view of \netheo{thm1} and \netheo{th2} it follows that $(Y_1,Y_2)$ is also an elliptically symmetric random vector with
some new random radius $R^*=RS$, and if $R$ is a Webullian-type rv satisfying \eqref{WeibG} with $p_2,L_2$ positive and some
$g_2$ being regularly varying at infinity, then $R^*$ has a distribution with scaling function $a(x)= p_2 L_2 x^{\theta-1}, \theta=p_2$ when $S$ is bounded, or
$$a(x)= \frac{B p_1p_2}{p_1+p_2}x^{\theta -1}, \quad \theta=\frac{p_1p_2}{p_1+p_2},$$
otherwise, where $B$ is given in \eqref{CA}. In view of Theorem 2.1 in Hashorva (2010) we obtain
$$ \bar \chi_{(Y_1,Y_2)}= 2 \Bigl(\frac{1 + \rho}{2}\Bigr)^{\theta/2} -1=: 2 \eta- 1,$$
where $\eta$ is the coefficient of tail dependence introduced in Ledford and Tawn (1996).

In \tb{a} particular case that $(Y_1,Y_2)$ has an elliptically generalized hyperbolic distribution with some underlying parameter $\rho \in (-1,1)$,
which was considered in Schlueter and Fischer (2012), then it follows easily that such model corresponds to a scaled elliptical random vector where
$S^2$ is some generalized inverse Gaussian rv, i.e., it has Weibullian tail behavior with $p_1=2$ and $L_1>0$.\\
 Since $R$ is also a Weibullian-type rv  with $p_2=2,L_2=1/2$, we get that $\theta=1$, and thus
$$ \bar \chi_{(Y_1,Y_2)}= \sqrt{ 2(1+ \rho)} -1=: 2 \eta- 1,$$
which is the result of Theorem 1 in Schlueter and Fischer (2012).\\
\tb{Note that} a correct form  of the joint asymptotic tail behavior
of $(Y_1,Y_2)$ (derived in Theorem 3 therein) can be directly obtained by applying Theorem 2.1 in Hashorva (2010).
}

\def\NN{\mathcal{N}}
\def\TT{\mathcal{T}}
\subsection{Supremum over random intervals for Gaussian processes with stationary increments}
The main result of Arendarczyk and D\c{e}bicki (2011) derives the exact asymptotics (as $u\to \IF$) of
$$ \pk{ \sup_{t \in [0, \mathcal{T}]} X(t)> u},$$
where $\{X(t),t\ge 0\}$ with $X(0)=0$ a.s.  is a mean-zero Gaussian process with stationary increments
and a.s. continuous sample paths being independent of $\mathcal{T}>0$,
which has tail asymptotics given by \eqref{poly}.
\tb{The following \tE{result} extends Theorem 3.1 in the aforementioned paper.}

\BT \label{th.sup}
\tc{
Let $\mathcal{T}$ be a nonnegative log-Weibullian rv
that satisfies (\ref{logW}) with some $L,p>0$}
and let
$\{X(t),t\ge 0\}$ be, an independent of $\mathcal{T}$,
centered Gaussian process with stationary increments
and continuously differentiable variance function
$\sigma^2(t)= Var(X(t))$. Suppose that $\sigma^2(\cdot)$ is convex,
regularly varying at infinity with index $\alpha \in (1,2]$. If further $\sigma^2(t) \le K t^\alpha$ holds for any $t>0$ and some positive constant $K>0$, then we have
\BQN
\pk{ \sup_{t \in [0, \mathcal{T}]} X(t)> u} \sim \pk{  \sigma(\mathcal{T}) \NN> u}, \quad u\to \IF,
\EQN
where $\mathcal{N}$ is \pE{an} $N(0,1)$ rv independent of $\TT$.
\ET
A combination of Theorem \ref{th2} with Theorem \ref{th.sup} leads to the following corollary.
\begin{corr}\label{corr.1}
Under the setup of \netheo{th.sup} suppose further that
$\sigma(t)\sim Ct^{\alpha/2}$ as $t\to\infty$ with $\alpha \in (1,2]$ and some constant $C>0$. \\
(a)
\tc{Then $\sigma(\mathcal{T})$ satisfies (\ref{logW}) with
$\widetilde{p}=\frac{2p}{\alpha}$ and $\widetilde{L}=\frac{L}{C^p}$
and
\BQN
\lim_{u\to\infty}
\frac{\log\left(\pk{ \sup_{t \in [0,\mathcal{T}]} X(t)> u} \right)}{u^{2\widetilde{p}/(\widetilde{p}+2)}}
=-
\widetilde{L}(\widetilde{L}\widetilde{p})^{\frac{-\widetilde{p}}{\widetilde{p}+2}} -
\frac{1}{2}(\widetilde{L}\widetilde{p})^{\frac{2}{\widetilde{p}+2}}=:-\widetilde{B}  .
\EQN
}
(b)\tc{
If
$\sigma(\mathcal{T})$ satisfies (\ref{giv})
with $\widetilde{p},\widetilde{L}$ and some regularly varying at infinity function $\widetilde{g}$,
then
\BQN
\pk{ \sup_{t \in [0,\mathcal{T}]} X(t)> u} \sim
(\widetilde{p}+2)^{-\frac{1}{2}} \widetilde{g}\left( (\widetilde{L}\widetilde{p})^{\frac{-1}{\widetilde{p}+2}}u^{\frac{2}{\widetilde{p}+2}}\right)
\exp\left(-\widetilde{B} u^{\frac{2\widetilde{p}}{\widetilde{p}+2}} \right), \quad u\to \IF.
\EQN
}
\end{corr}

\begin{rem}
Clearly, if we specify in the assumptions  of Theorem \ref{th.sup}
that $\sigma^2(x)= C x^\alpha$ (i.e., $X$ is a fractional Brownian motion
with Hurst index $\alpha/2$)
and $\mathcal{T}$ is Weibullian, then both $\sigma(\TT)$ and $\NN$ are Weibullian-type rvs,
and thus assumptions of Corollary \ref{corr.1} (b) are satisfied. Hence
Corollary \ref{corr.1} is an extension of Theorem 4.1 in Arendarczyk and D\c{e}bicki (2011).
\end{rem}

\section{Proofs}
It  is well-known that for some rv $U$ which has distribution with right endpoint equal to infinity the assumption $U \in GMDA(a)$  implies that the tail of $U$ is rapidly varying at infinity, i.e.,
$$ \lim_{u \to \IF} \frac{\pk{U> \lambda u}}{\pk{U> u}}= 0$$
holds for any $\lambda> 1$. First we present a result on the random scaling of rvs with rapidly varying tails which is of some interest on its own. 

\BEL \label{XXH} Let $S, Y,Y^*$ be three independent rvs. Suppose that $S\ge 0$ has distribution $G$ with right endpoint equal to 1. If further $Y$ has a rapidly varying tail and $\pk{Y> u} \sim L(u) \pk{Y^*> u}$ as $u\to \IF$ for some slowly varying function $L(\cdot)$, then
\BQN
\pk{S Y > u} \sim \pk{SY>u, S > w} \sim  L(u)\pk{SY^*>u}
\EQN
holds for any $w\in (0,1)$.
\EEL

\prooflem{XXH}
By the independence of $S$ and $Y$ for any $u>0$ and $w\in (0,1)$, we have
\BQNY
 \pk{SY> u}
 &=& \pk{S Y> u, S \le w} + \pk{S Y> u, S > w}\\
 &\le & \pk{ Y> u/w}+ \pk{S Y> u, S > w}. 
 \EQNY
The assumption that  $Y$ has  a rapidly varying tail implies for any $t\in (w,1)$
\BQNY
\frac{\pk{Y> u/w}}{\pk{S Y > u, Y > w}} &\le &\frac{ \pk{Y> u/w} }{ \pk{S Y > u, S > t}} \le  \frac{\pk{Y> u/w}}{ \pk{Y > u/t} \pk{S > t}} \to 0, \quad u\to \IF
\EQNY
hence for any $w\in (0,1)$
\BQNY\label{simp}
 \pk{SY> u} \sim \int_w^1 \pk{Y> u/s} \, d G(s), \quad u\to \IF.
 \EQNY
By the uniform convergence theorem for regularly varying functions (e.g., Embrechts et al.\ (1997))
 \BQNY
 \int_w^1 \pk{Y> u/s} \, d G(s) \sim L(u)\int_w^1 \pk{Y^*> u/s} \, d G(s), \quad u\to \IF.
 \EQNY
 The assumption $\pk{Y> u} \sim L(u) \pk{Y^*> u}$ as $u\to \IF$ yields that  $Y^*$ has also a rapidly varying tail at infinity. Hence  in view of the above arguments and the fact that  $S$ and $Y^*$ are independent, we have that
 $$ \int_w^1 \pk{Y^*> u/s} \, d G(s) \sim \pk{SY^*> u}, \quad u\to \IF$$
 establishing the proof. \QED

\prooftheo{thm1} The assumption $\pE{Y_1}\in GMDA(a)$ implies that the convergence
\BQN\label{drit}
\frac{\pk{Y_1>u+ x a(u)}}{\pk{Y_1> u}} \to \exp(-x), \quad u\to \IF
\EQN
holds uniformly for $x$ \pE{on} compact sets of $\R$.
\COM{Since $Y_1$ and $Y_2$ are independent for any $w\in (0,1)$, we have
\BQNY
 \pk{Y_1Y_2> u}= \int_0^1 \pk{Y_1> u/s} \, d G(s)&=& \pk{Y_1 Y_2> u, Y_2 \le w} + \pk{Y_1 Y_2> u, Y_2 > w}\\
 &\le & \pk{ Y_1> u/w}+ \int_w^1\pk{Y_1> \pE{u/s}} \, d G(s),
 \EQNY
 with $G$ the distribution  of $Y_2$.
Hence, since \aa{\eqref{drit} yields}
$$ \limit \frac{ \pk{Y_1> u/w}}{\pk{Y_1> u/t}}=0,$$
we get for some $t\in (w,1)$
\BQNY
\frac{\pk{Y_1> u/w}}{\pk{Y_1 Y_2 > u, Y_2 > w}} &\le &\frac{ \pk{Y_1> u/w} }{ \pk{Y_1 Y_2 > u, Y_2 > t}} \\
&\le & \frac{\pk{Y_1> u/w}}{ \pk{Y_1 > u/t} \pk{Y_2 > t}} \to 0, \quad u\to \IF,
\EQNY
 hence for any $w\in (0,1)$
\BQN\label{simp}
 \pk{Y_1Y_2> u} \sim \int_w^1 \pk{Y_1> u/s} \, d G(s), \quad u\to \IF.
 \EQN
 }
Since $Y_1$ has a rapidly varying tail at infinity, then  by \nelem{XXH} for any \aa{fixed} $z\ge 0$
 and $w\in (0,1)$
 \BQNY
\pk{Y_1Y_2 > u+ a(u)z } \sim \int_{w}^1 \pk{Y_1>(u+ z a(u))/s}\, d G(s), \quad u\to \IF
\EQNY
holds with $G$ the distribution  of $Y_2$. By the uniform convergence theorem for regularly varying functions
$$ \limit \frac{a(ux)}{a(u)}= x^{-\tau}$$
holds uniformly
\tc{for $x\in[w,1]$,}
with $w\in (0,1)$ some arbitrary constant.  \aa{Hence}
$$ z_{u,s}:= \frac{z}{s} \frac{a(u)}{a(u/s)} \to \frac{z}{s^{1+\tau}}, \quad u\to \IF$$
uniformly for $s\in [w,1]$, \aa{and thus}
$$ \frac{\pk{Y_1>u/s+ z a(u)/s}}{\pk{Y_1>u/s}} = \frac{\pk{Y_1>u/s+  a(u/s) z_{u,s}}}{\pk{Y_1>u/s}} \to \exp(- z/ s^{1+\tau}), \quad u\to \IF$$
uniformly for $s\in [w ,1]$. For any $\ve>0$ we can find $w \in (0,1)$ such that for all $s\in [w,1]$
$$ (1- \ve) \exp(-z)\le \exp(- z/ s^{1+\tau}) < (1+ \ve) \exp(-z)$$
implying that as $u\to \IF$
$$ \pk{Y_1Y_2> u+ a(u)z} \sim \int_{w}^1 \pk{Y_1> u/s+ a(u/s) z_{u,s}}\, d  G(s) \sim \exp(-z) \pk{Y_1Y_2> u}.$$
Hence $Y_1Y_2\in GMDA(a)$, and thus the proof is complete.
\QED

\prooftheo{th2} {\it Ad.(a).} 
Since for any $u>0$
\begin{eqnarray*}
\pk{Y_1Y_2>u}
\ge
\pk{Y_1>\left(\frac{p_2L_2}{p_1L_1}\right)^{1/(p_1+p_2)}u^{p_2/(p_1+p_2)}}
\pk{Y_2>\left(\frac{p_1L_1}{p_2L_2}\right)^{1/(p_1+p_2)}u^{p_1/(p_1+p_2)}},
\end{eqnarray*}
then we immediately get
\[
\liminf_{u\to\infty}
\frac{\log(\pk{Y_1Y_2>u})}{u^{p_1p_2/(p_1+p_2)}}
\ge
-\left(
   L_1\left(\frac{p_2L_2}{p_1L_1}\right)^{\frac{p_1}{p_1+p_2}}
  +
   L_2\left(\frac{p_1L_1}{p_2L_2}\right)^{\frac{p_2}{p_1+p_2}}
   \right)=:-\pE{B}.
\]
Next, we 
have
\begin{eqnarray*}
\pk{Y_1Y_2>u}
&\le&
\sum_{k=[u^{p_2/(2(p_1+p_2))}]}^{[u^{\tc{(p_2+p_1/2)}/(p_1+p_2)}]}
\pk{Y_1\in [k,k+1), Y_1Y_2>u}\\
&&+
\pk{Y_1<[u^{p_2/(2(p_1+p_2))}], Y_1Y_2>u}
+
\pk{Y_1>[u^{\tc{(p_2+p_1/2)}/(p_1+p_2)}], Y_1Y_2>u}\\
&=&
\Sigma+P_1+P_2.
\end{eqnarray*}
Now observe that, as $u\to\infty$
\begin{eqnarray}
\log\left(P_1\right)
\le
\log\left(\pk{Y_2>u^{1-p_2/(2(p_1+p_2))}}\right)
\sim -L_2u^{(p_1+p_2/2)p_2/(p_1+p_2)}\label{l.neg1}
\end{eqnarray}
and
\begin{eqnarray}
\log\left(P_2\right)
\le
\log\left(\pk{Y_1>[u^{(p_2+p_1/2)/(p_1+p_2)}]}\right)
\sim -\tc{L_1 u^{(p_2+p_1/2)p_1/(p_1+p_2)}}\label{l.neg2}.
\end{eqnarray}
Moreover, for each $\varepsilon>0$, sufficiently large $u$ and
$k\in\left[[u^{p_2/(2(p_1+p_2))}],[\tc{u^{(p_2+p_1/2)/(p_1+p_2)}}]\right]$
\begin{eqnarray}
\log\left(\pk{Y_1\in [k,k+1), Y_1Y_2>u}\right)
&\le&
\log\left(\pk{Y_1\ge k, Y_2>u/(k+1)}\right)\nonumber\\
&\le&
\tc{-(1-\varepsilon)\left(L_1k^{p_1}+L_2(u/k)^{p_2}\right)\nonumber}\\
&\le&
- (1-\varepsilon)\pE{B}
   u^{p_1p_2/(p_1+p_2)}, \label{explan1}
\end{eqnarray}
\tc{
where (\ref{explan1}) follows from the fact that
$f(x)=L_1x^{p_1} +L_2\left(\frac{u}{x}\right)^{p_2}$
attains its minimum $f(x_u)=B   u^{p_1p_2/(p_1+p_2)}$ at $x_u=\left( \frac{p_2L_2}{p_1L_1} \right)^{1/(p_1+p_2)}u^{p_2/(p_1+p_2)}$
and for any $\delta\in (0,1)$ and all $u$ large $k/(k+1)> 1- \delta$.}  
Thus, using \aa{the fact} that $\Sigma$ consists of a polynomial (with respect to $u$)
number of elements,
we have that
\begin{eqnarray}
\limsup_{u\to\infty}
\frac{\log(\Sigma)}{u^{p_1p_2/(p_1+p_2)}}\le -\pE{B}.\label{non.neg}
\end{eqnarray}
The combination of (\ref{l.neg1}), (\ref{l.neg2}) with (\ref{non.neg})
completes the proof of \aa{statement}  (a).
\\
{\it Ad. (b).}  Suppose without loss of generality that $L_1=L_2=1$. With the same arguments as in the proof of Lemma 3.2 in Hashorva and Ji (2013),
if $Y_1^*$ and $Y_2^*$ are two positive independent rvs tail equivalent to $Y_1$ and $Y_2$, respectively,  then
$$ \pk{Y_1Y_2> u} \sim \pk{Y_1^*Y_2^*>u}, \quad u\to \IF.$$
We define next $Y_i^*=S_iZ_i$ where $S_i$ has distribution $G_i,i=1,2$ with right endpoint equal to 1, and $Z_1,Z_2$ are independent of $S_1,S_2$.
 Let $\alpha_1^*$ and $\alpha_2^*$ be the index of the  regular variation of $g_1$ and $g_2$, respectively.  Let $\alpha_i> \alpha_i^*,i=1,2$ be two arbitrary constants.
 The functions $\tilde g_i(x)= g_i(x) x^{-\alpha_i}$ are regularly varying at infinity with index $\alpha_i^*- \alpha_i < 0$.
 Hence, we can assume without loss of generality, that
 $$ \pk{S_i > 1 -  a_i(u)/u} = \frac{1}{\Gamma(\alpha_i- \alpha_i^*+1)} \tilde g_i(u), \quad u\to \IF,$$
 where $a_i(u)= u^{1- p_i}, i=1,2,u>0.$  
 In view of  Example 1 in Hashorva (2012) (see also Theorem 3.1 in Hashorva et al.\ (2010))  for $i=1,2$ we obtain
 $$ \pk{S_i Z_i> u}   \sim \pk{Y_i> u}, \quad u \to \IF,$$
where $S_i,Z_i, i=1,2$ are independent and positive rvs, and
 $$ \pk{Z_i> u}\sim u^{\alpha_i} \exp(-  u_i^{p_i}), \quad u\to \IF.$$
Consequently, as $u\to \IF$
\BQNY
 \pk{Y_1Y_2 > u} &\sim &\pk{S_1Z_1 S_2 Z_2> u} \sim \pk{ U W > u},
 \EQNY
 where $U=S_1S_2$ and $W=Z_1Z_2$. The tail asymptotics of $U$ follows by
 a direct application of Theorem 2.1 in Farkas and Hashorva (2013) whereas the tail asymptotics of
 $W$ follows from \eqref{DE}. Hence, the tail asymptotics of $UW$ follows by applying again the result of the aforementioned example, and thus the proof is complete.
 \QED

\prooftheo{th:last} The proof follows \tb{by} the same arguments as the proof of Theorem 5.1 in Hashorva (2013).
When $S$ is a bounded rv, then in view of Remark \ref{RemA}  we have that $d_n= (1+o(1))\sqrt{2  \log  n}$ and since the
scaling function $a(x)=1/x$, then $c_n=1/d_n$ follows. For the case $S$ has a Weibullian tail behavior, the relation between $c_n$ and $d_n$ can be established using the same idea as in the proof of the aforementioned theorem. \QED

\prooftheo{th.sup} For chosen constants $\gamma_1=2/(\alpha+2p)$, $\gamma_2=4/(2\alpha+p)$
and
$\delta=\delta(u)=2\frac{\sigma^3(u)}{\sigma'(u)}u^{-2}\log^2(u)$ we have
\begin{eqnarray*}
\pk{\sup_{t\in[0,\mathcal{T}]} X(t)>u}
&\le&
\int_0^{u^{\gamma_1}}
 \pk{\sup_{t\in[0,s]} X(t)>u}dF_{\mathcal{T}}(s)
 +
 \int_{u^{\gamma_1}}^{u^{\gamma_2}}
 \pk{\sup_{t\in[0,s-\delta]} X(t)>u}dF_{\mathcal{T}}(s)\\
 &&+
 \int_{u^{\gamma_1}}^{u^{\gamma_2}}
 \pk{\sup_{t\in[s-\delta,s]} X(t)>u}dF_{\mathcal{T}}(s)
 +
 \int_{u^{\gamma_2}}^\infty
 \pk{\sup_{t\in[0,s]} X(t)>u}dF_{\mathcal{T}}(s)\\
&=&
I_1+I_2+I_3+I_4.
\end{eqnarray*}
\tE{As in} \COM{Following line-by-line the same argument as in} the proof of
Theorem 3.1 in Arendarczyk and D\c ebicki (2011),
we conclude that
\[
I_1+I_2=o(\pk{X(\mathcal{T})>u})
\]
as $u\to\infty$ and for each $\varepsilon>0$ and $u$ large enough
\[
I_3\le (1+\varepsilon)\pk{X(\mathcal{T})>u}=(1+\varepsilon)\pk{\sigma(\TT)\mathcal{N}>u},
\]
where  $\mathcal{N}$ is \pE{an} $N(0,1)$ rv independent of $\TT$. Thus it suffices to show that
\begin{eqnarray}
I_4=o(\pk{X(\mathcal{T})>u})\label{sm.I4}
\end{eqnarray}

as $u\to\infty$.
Indeed, since for all large $u$ we have
$I_4\le \pk{\mathcal{T}>u^{\gamma_2}},$
  then
\[
\limsup_{u\to\infty}\frac{\log(I_4)}{u^{4p/(2\alpha+p)}}\le -L.
\]
On the other hand,
for each $\varepsilon\in(0,\alpha/2)$
and sufficiently large $u$, the assumption that $\sigma(\cdot)$
is regularly varying at $\infty$ with index $\alpha/2$ implies
\[
\pk{\sigma(\mathcal{T})>u}\ge \pk{\mathcal{T}^{\alpha/2-\varepsilon}>u}.
\]
Hence, \tE{for some $K>0$} by statement (a) of Theorem \ref{th2}
\[
\liminf_{u\to\infty}
\frac{\log(\pk{X(\mathcal{T})>u})}
      {u^{2p/(\tc{p}+\alpha-2\varepsilon)}}
\ge
\liminf_{u\to\infty}
\frac{\log(\pk{\mathcal{T}^{\alpha/2-\varepsilon}\mathcal{N}>u})}
      {u^{2p/(\tc{p}+\alpha-2\varepsilon)}}\ge -K.
\]
Consequently, since for sufficiently small $\varepsilon>0$,
we have $2p/(\tc{p}+\alpha-2\varepsilon)<4p/(p+2\alpha)$, then
(\ref{sm.I4}) holds.
\QED

\COM{\proofprop{infsum} 
By the assumptions we have $$S_m:= \sum_{i=1}^m \lambda_i R_i X_i Y_i= R_1 \sum_{i=1}^m X_iY_i \equaldis R_1 X_1 \sqrt{ \sum_{i=1}^m Y_i^2},$$
where $\equaldis$ stands for equality of distributions. Since $X_1$ is symmetric about 0, and
$\abs{X_1}$ and $\sqrt{ \sum_{i=1}^2 Y_i^2}$ are both Weibullian-type rvs, by the independence of $R_1, X_i, Y_i,i\le m$
we have by \netheo{th2} that $\abs{X_1} \sqrt{ \sum_{i=1}^m Y_i^2}$ is also a Weibullian-type rv (see also Lemma 3.4 in Hashorva et al.\ (2012)). Furthermore, since $R_1$ is bounded, by \netheo{thm1}
$\abs{S_m}$ is also Weibullian-type and we obtain
$$ S_m \in GMDA(a), \quad a(x)=1/x.$$
Now, the claim follows by using the same arguments as those in the proof of Theorem 2.1 in the aforementioned paper. \QED
}

\proofkorr{corr.1} The proof
boils down to checking, that
for both  cases (a) and (b) the conditions imposed  on $\sigma(\cdot)$ imply that
$\TT$ satisfies the assumptions of Theorem \ref{th.sup};  therefore we omit the details. \QED.

{\bf Acknowledgement}:
K. D\c{e}bicki was supported by NCN Grant No 2011/01/B/ST1/01521 (2011-2013).
The authors acknowledge partial support by the project RARE -318984, an FP7 IRSES Fellowship and the Swiss National Science Foundation Grant 200021-140633/1.

\bibliographystyle{plain}

\begin{thebibliography}{99}
\bibitem{} Arendarczyk, M., D\c{e}bicki, K. (2011) Asymptotics of supremum distribution of
a Gaussian process over a Weibullian time.  \emph{Bernoulli}, {\bf 17}, 194--210.

\bibitem{} Arendarczyk, M., D\c{e}bicki, K. (2012) Exact asymptotics of supremum of a stationary Gaussian process over a random interval.
 \emph{Stat. Probab. Letters},  {\bf 82}, 645--652.


\bibitem{BERMANc}
\cE{Berman,  M.S. (1992) {\it Sojourns and Extremes of Stochastic Processes}, Wadsworth \& Brooks/ Cole, Boston.}





\bibitem{CLSG}
Cline, D.B.H.,  Samorodnitsky, G. (1994) Subexponentiality of the product of independent random variables.
\emph{Stoch. Proc. Appl.} {\bf 49}, 1, 75--98.


\bibitem{}
\tE{Charpentier, A., Segers, J. (2007) Lower tail dependence for
Archimedean copulas: characterizations and pitfalls. \emph{Insurance
Math. Econom.} {\bf 40}, 525--532.}

\bibitem{}
\tE{Charpentier, A., Segers, J. (2009) Tails of multivariate
Archimedean copulas. \emph{J. Multivariate Anal.} {\bf 100},
1521--1537.}

\tE{
\bibitem{}
D'Auria, B.,  Resnick, S.I. (2006) Data network models of
burstiness. \emph{Adv. in Appl. Probab.} {\bf 38}, 373--404.
}

\tE{\bibitem{}
D'Auria, B.,  Resnick, S.I. (2008) The influence of dependence
on data network models. \emph{Adv. in Appl. Probab.} {\bf 40},
60--94.
}



\bibitem{} D\c{e}bicki, K., Hashorva, E., Ji, L. (2013)
Tail asymptotics of supremum of certain Gaussian processes over threshold dependent random intervals. Extremes, in press. 

\bibitem{} D\c{e}bicki, K., van Uitert, M. (2006)  Large buffer asymptotics for generalized processor sharing queues
with Gaussian inputs. \emph{Queueing Syst},  {\bf 54}, 111--120.

\bibitem{} D\c{e}bicki, K., Zwart, A.P., Borst, S.C. (2004) The supremum of a Gaussian process
over a random interval.  \emph{Stat. Probab. Letters}, {\bf 68}, 221--234.


\bibitem{}
Embrechts, P., Kl\"uppelberg, C., Mikosch, T. (1997) {\it Modeling Extremal Events For Finance And Insurance.}
Berlin, Springer.


\bibitem{farkas}
Farkas, J., Hashorva, E. (2013) Tail approximation for reinsurance portfolios of Gaussian-like risks.
\emph{Scand. Actuar. J.}, in press, DOI 10.1080/03461238.2013.825639.



\bibitem{}
Hashorva, E. (2013) Minima and maxima of elliptical triangular arrays and spherical processes. \emph{Bernoulli},
{\bf 19}, 886--904.

\bibitem{}
Hashorva, E., Ji, L. (2013) Asymptotics of the finite-time ruin probability for the Sparre Andersen risk model perturbed by an inflated stationary chi-process. \emph{Communications in Statistics - Theory and Methods}, {\bf 43}, 2540-–2548.


\bibitem{}
 Hashorva,  E.,   Weng, Z. (2013) Tail asymptotic of Weibull-type risks. \emph{Statistics}, in press
 DOI:10.1080/02331888.2013.800520.

\bibitem{}
 Hashorva,  E.  (2012) Exact tail asymptotics in bivariate scale mixture models. \emph{Extremes}, {\bf 15},  109--128.


\bibitem{}
 Hashorva,  E., Pakes, A.G., and Tang, Q. (2010) Asymptotics of random contractions. 
  \emph{Insurance: Mathematics and Economics,} {\bf 47},  405--414.




\bibitem{} Jessen, A.H.,  Mikosch, T. (2006) Regularly varying functions. \emph{Institut Math\'{e}matique.
Publications. Nouvelle S\'{e}rie.} {\bf 80}, 171--192.


\bibitem{}
Kabluchko, Z. (2011) Extremes of independent Gaussian processes. \emph{Extremes}, {\bf 14}, 285--310.


\bibitem{} Maulik, K., Resnick, S. (2004)
Characterizations and examples of hidden regular variation. {\it Extremes}, {\bf 7}, 31--67.

\bibitem{}
Nadarajah, S. (2005) Sums, products and ratios of generalized beta variables. \emph{Stat. Papers.} {\bf 47}, 69--90.

\bibitem{}
Pakes, A.G., Navarro, J. (2007)
Distributional characterizations through scaling relations. \emph{Aust. \& N.Z. J. Stat.} {\bf 49}, 115--135.


\bibitem{resnick1}
Resnick, S.I. (2007) {\it Heavy-Tail Phenomena: Probabilistic and
Statistical Modeling}. Springer, New York.

\bibitem{resnick2}
Resnick, S.I. (1987) {\it Extreme Values, Regular Variation and Point Processes.} \emph{Springer, New York}.

\bibitem{} Schlueter, S., Fischer, M. (2012) The weak tail dependence coefficient of the elliptical generalized hyperbolic
distriubtion. \emph{Extremes}. {\bf 15}, 159--174.


\bibitem{} Tan, Z., Hashorva, E.  (2013)
Exact asymptotics and limit theorems for supremum of stationary chi-processes over a random interval.
\emph{Stoch. Proc. Appl.} {\bf 123}, 2983--2998.


\bibitem{} Yang, Y., Wang, Y. (2013) Tail behavior of the product of two dependent random variables with applications to risk theory.  \emph{Extremes}, {\bf 16}, 55--74.


\end{thebibliography}

\end{document}